\documentclass[11pt]{article}
\usepackage{amssymb,amsfonts,amsmath,latexsym,epsf,tikz,url}

\newtheorem{theorem}{Theorem}[section]

\newtheorem{corollary}[theorem]{Corollary}
\newtheorem{lemma}[theorem]{Lemma}

\newcommand{\proof}{\noindent{\bf Proof.\ }}
\newcommand{\qed}{\hfill $\square$\medskip}

\begin{document}

\title{Distinguishing number and distinguishing index of  strong product of two graphs}

\author{ Samaneh Soltani \and 
Saeid Alikhani  $^{}$\footnote{Corresponding author}
}

\date{\today}

\maketitle

\begin{center}
Department of Mathematics, Yazd University, 89195-741, Yazd, Iran\\
{\tt  s.soltani1979@gmail.com, alikhani@yazd.ac.ir}
\end{center}

%%%%%%%%%%%%%%ABSTRACT%%%%%%%%%%%%%%%%%%%%%%%%%%%%%%%%%%%%%%%%%%%%%%%%%%%%%%%%%%%%

\begin{abstract}
The distinguishing number (index) $D(G)$ ($D'(G)$)  of a graph $G$ is the least integer $d$
such that $G$ has an vertex labeling (edge labeling)  with $d$ labels  that is preserved only by a trivial automorphism. The strong product $G\boxtimes H$ of two graphs $G$ and $H$ is the graph with vertex set $V (G)\times V (H)$ and edge set $\{\{(x_1, x_2), (y_1, y_2)\} | x_iy_i \in E(G_i) ~{\rm or}~ x_i = y_i ~{\rm for~ each}~ 1 \leq i \leq 2.\}$.
In this paper we study the distinguishing number and the distinguishing index of strong product of two graphs.  We prove that for every $k \geq 2$, the k-th strong power of a connected $S$-thin
graph $G$ has distinguishing index equal 2.

\end{abstract}


\medskip
\noindent{\bf AMS Subj.\ Class.:} 05C15, 05C60. 

%%%%%%%%%%%%%%%%%%%%%%%%%%%%%%%%%%%%%%%%%%%%%%%%%%%%%%%%%%%%%%%%%%%%%%%%%%%%%%%%%
%%%%%%%%%%%%%%%%%%%%%%%%%%%%%%%%%%%%%%%%%%%%%%%%%%%%%%%%%%%%%%%%%%%%%%%%%%%%%%%%%
\section{Introduction and definitions}
Let $G=(V,E)$ be a simple graph of order $n\geq 2$. We use the the following notations: The set of vertices adjacent in $G$ to a vertex of a vertex subset $W\subseteq  V$ is the \textit{open neighborhood}  $N_G(W )$ of $W$. Also $N_G(W )\cup W$ is called a \textit{closed neighborhood}   of $W$ and denoted by $N_G[W]$. ${\rm Aut}(G)$ denotes the automorphism group of $G$.  A  \textit{subgraph} of a graph $G$ is a graph $H$ such that 
$V(H) \subseteq V(G)$ and  $E(H) \subseteq E(G)$. 
If $V(H) = V(G)$, we call $H$ a \textit{spanning subgraph} of $G$. Any spanning 
subgraph of $G$ can be obtained by deleting some of the edges from $G$.

A labeling of $G$, $\phi : V \rightarrow \{1, 2, \ldots , r\}$, is said to be \textit{$r$-distinguishing}, 
if no non-trivial  automorphism of $G$ preserves all of the vertex labels.
The point of the labels on the vertices is to destroy the symmetries of the
graph, that is, to make the automorphism group of the labeled graph trivial.
Formally, $\phi$ is $r$-distinguishing if for every non-trivial $\sigma \in {\rm Aut}(G)$, there
exists $x$ in $V$ such that $\phi(x) \neq \phi(\sigma(x))$. The \textit{distinguishing number} of a graph $G$ is defined  by
\begin{equation*}
D(G) = {\rm min}\{r \vert ~ G ~\textsl{\rm{has a labeling that is $r$-distinguishing}}\}.
\end{equation*} 

This number has defined by Albertson and Collins \cite{Albert}. Similar to this definition, Kalinowski and Pil\'sniak \cite{Kali1} have defined the \textit{distinguishing index}  $D'(G)$ of $G$ which is  the least integer $d$
such that $G$ has an edge colouring   with $d$ colours that is preserved only by a trivial
automorphism. If a graph has no nontrivial automorphisms, its distinguishing number is  $1$. In other words, $D(G) = 1$ for the asymmetric graphs.
The other extreme, $D(G) = \vert V(G) \vert$, occurs if and only if $G = K_n$. The distinguishing index of some examples of graphs was exhibited in \cite{Kali1}. For 
instance, $D(P_n) = D'(P_n)=2$ for every $n\geq 3$, and 
$D(C_n) = D'(C_n)=3$ for $n =3,4,5$,  $D(C_n) = D'(C_n)=2$ for $n \geq 6$. % It is easy to see that the value $|D(G)-D'(G)|$ can be large. For example $D'(K_{p,p})=2$ and $D(K_{p,p})=p+1$, for $p\geq 4$.  
%A graph and its complement, always have the same automorphism group while their graph structure usually differs, hence $D(G) = D(\overline{G})$ for every simple graph $G$.
The distinguishing number and the distinguishing index of some graph products has been studied in literature (see \cite{Imrich,Imrich and  Klavzar}). 
The \textit{Cartesian product} of graphs $G$ and $H$ is a graph, denoted $G\Box H$, whose vertex
set is $V (G) \times V (H)$. Two vertices $(g, h)$ and $(g', h')$ are adjacent if either $g = g'$ and
$hh' \in E(H)$, or $gg' \in E(G)$ and $h = h'$. 
The \textit{direct product}  of $G$ and $H$ is the graph, denoted as $G \times H$, whose vertex set is
$V (G) \times V (H)$, and for which vertices $(g, h)$ and $(g_0, h_0)$ are adjacent precisely if $gg_0 \in E(G)$
and $hh_0 \in E(H)$. The \textit{strong product} of $G$ and $H$ is the graph denoted as $G \boxtimes H$, and defined by
\begin{equation*}
V (G \boxtimes H) = \{(g, h) | g \in V (G) ~{\rm and}~ h \in V (H)\},~E(G \boxtimes H) = E(G\Box H) \cup E(G \times  H).
\end{equation*}
 Note that $G\Box H$ and $G\times H$ are  spanning subgraphs of $G \boxtimes H$. Note that $K_n \boxtimes K_m =K_{nm}$, and $G\boxtimes K_1 =G$.

Let $*$ represent either the Cartesian or the strong product operation, and
consider a product $G_1 * G_2 * \cdots *  G_k$. For any index $1 \leq i \leq k$, there is a projection map $p_i :  G_1 * G_2 * \cdots *  G_k \rightarrow G_i$ defined as $p_i(x_1, x_2, \ldots  , x_k) = x_i$. We call $x_i$ the $i$th coordinate
of the vertex $(x_1, x_2, \ldots  , x_k)$. Given a vertex $a=(a_1, a_2, \ldots  , a_k)$ of the product $G = G_1 * G_2 * \cdots *  G_k$, the
$G_i$-layer through $a$ is the induced subgraph
\begin{equation*}
G^a_i = \{x \in V (G) | p_j(x) = a_j ~{\rm for}~ j \neq i\}= \{(a_1, a_2, \ldots , x_i, \ldots , a_k) | x_i \in V (G_i)\}.
\end{equation*}
Note that $G^a_i = G^b_i$ if and only if $p_j(a) = p_j(b)$ for each index $j \neq i$.

A graph is prime with respect to a given graph product if it is nontrivial and cannot be
represented as the product of two nontrivial graphs. For the Cartesian and strong product, this means that a nontrivial graph $G$ is prime if $G = G_1 \Box G_2$ ($G = G_1 \boxtimes G_2$ ) implies that $G_1$ or $G_2$ is $K_1$.
 We need knowledge of the structure of the automorphism group of the Cartesian product, which was determined by Imrich \cite{Imrich1969}, and independently by Miller \cite{Miller}.

\begin{theorem}{\rm \cite{Imrich1969,Miller}}\label{autoCartesian} Suppose $\psi$  is an automorphism of a connected graph $G$ with prime
	factor decomposition $G = G_1\Box G_2\Box \ldots \Box G_r$. Then there is a permutation $\pi$ of the set
	$\{1, 2, \ldots ,  r\}$ and there are isomorphisms  $\psi_i : G_{\pi(i)} \rightarrow G_i$, $i = 1,\ldots , r$, such that
	$$\psi(x_1, x_2 , \ldots , x_r) = ( \psi_1(x_{\pi(1)}),  \psi_2(x_{\pi(2)}), \ldots ,  \psi_r(x_{\pi(r)})).$$ 
\end{theorem}

 We say vertices $x, y$ of a graph $G$ are in the \textit{ relation $S$}, written
$xSy$, provided that $N[x] = N[y]$. We recall that
graphs with no pairs of vertices with the same closed neighborhoods are called $S$-thin. The relationship between the automorphism group of a connected $S$-thin graph and the groups of its prime factors with respect to the strong product is the same as that in the case of the Cartesian product.
\begin{theorem}{\rm \cite{Sandi}}\label{autstrong} Suppose $\varphi$ is an automorphism of a connected  $S$-thin graph $G$ that has a prime factorization $G = G_1 \boxtimes G_2 \boxtimes \ldots \boxtimes  G_k$. Then there exists a permutation 	$\pi$ of $\{1, 2, \ldots , k\}$, together with isomorphisms $\varphi_i : G_{\pi(i)}\rightarrow G_i$, such that $$\varphi(x_1, x_2, \ldots , x_k) = (\varphi_1(x_{\pi(1)}), \varphi_2(x_{\pi(2)}), \ldots , \varphi_k(x_{\pi(k)})).$$
	 \end{theorem}
In other words,  there exists a $G_i$ to every $G_j$, $1 \leq i, j \leq k$, such
that $\varphi$ maps every $G_j$-layer into a $G_i$-layer.

\medskip

In the next section we study the distinguishing number of strong product of two graphs. In section 3, we show that the distinguishing index of strong product of two simple connected, $S$-thin prime graphs cannot be more than the distinguishing index of their direct product. As a consequence we prove that all powers of a  $S$-thin graph $G$ with respect to the strong product distinguished by exactly two edge labels.

\section{Distinguishing number of strong product of two graphs}
We begin this section with a general upper bound for the strong product of two simple connected graphs.
\begin{theorem}\label{generalbound}
If $G$ and $H$ are  two simple connected graphs, then 
\begin{equation*}
D(G\Box H)\leq D(G\boxtimes H)\leq  {\rm min}\left\{D(G)|V(H)|, |V(G)|D(H)\right\}.
\end{equation*}
\end{theorem}
\proof For the left inequality it is sufficient to know that  $G\Box H$ is a spanning subgraph of $G\boxtimes H$ such that ${\rm Aut}(G\Box H)\subseteq {\rm Aut}(G\boxtimes H)$. For the right inequality, we present two distinguishing vertex labelings of $G\boxtimes H$ with $D(G)|V(H)|$ and $|V(G)|D(H)$ labels,  respectively. For the first labeling, let $\phi$ be a $D(G)$-distinguishing labeling of the vertices of $G$. We assign the vertex $(g_i ,h_j)$ the label $\phi(g_i)+(j-1)D(G)$ where $1\leq i \leq |V(G)|$ and $1\leq j \leq |V(H)|$. In fact, we label each $G$-layer with $D(G)$ different labels such that the sets of $D(G)$ labels  which have been assigned  to each $G$-layer does not have intersection. Hence if $f$ is an automorphism of $G\boxtimes H$, then $f$ maps each $G$-layer to itself, and so $f$ is a automorphism on $G$-layers. Now since every $G$-layer is labeled distinguishingly, so $f$ is the identity automorphism. Therefore we have a $D(G)|V(H)|$-distinguishing labeling of $G\boxtimes H$. Repeating the above argument with $H$-layers instead of $G$-layers, yields a $|V(G)|D(H)$-distinguishing labeling.\qed

The right bound of Theorem \ref{generalbound} is sharp for strong product of complete graphs.

Since the automorphism group of a strong product of connected, $S$-thin prime graphs is the same as automorphism group of the Cartesian product of them, so the following theorem follows immediately:

\begin{theorem}\label{disnumstroncartes}
If $G$ and $H$ are two simple connected, $S$-thin prime graphs, then $D(G\boxtimes H)=D(G\Box H)$.
\end{theorem}

Since the path graph $P_n$ ($n\geq3$), and the cycle graph $C_m$ ($m\geq5$) are connected, $S$-thin prime graphs, then by Theorem \ref{disnumstroncartes} we have $D(P_n \boxtimes P_q)= D(P_n \boxtimes C_m)=D(C_m \boxtimes C_p)=2$ for any $q,n\geq3$ and $m,p\geq5$. (see \cite{Imrich and  Klavzar} for the distinguishing number of Cartesian product of these graphs).

Imrich and  Kla\v{v}zar in \cite{Imrich and  Klavzar} proved that the distinguishing number of  $k$-th power  with respect to the strong product of a  connected, $S$-thin graph  is two.

\begin{theorem}{\rm \cite{Imrich and  Klavzar}}
	Let $ G$ be a  connected, $S$-thin graph	and $\boxtimes G^k$ the $k$-th power of $G$ with respect to the strong product. Then $D(\boxtimes G^k) = 2$ for $k \geq 2$.
\end{theorem}

%\begin{theorem}
%Let  $G$ and $H$ be two nonisomorphic simple connected, $S$-thin prime graphs of orders $n$ and $m$, respectively. If $m\leq 2^n$ and $D(G)\geq2$, then $D(G\boxtimes H)\leq D(G)$.
%\end{theorem}
%\proof We have at most $2^n$ distinct $n$-ary of labels with two labels 1 and 2, which we denote them by $a_i=(a_{i1},\ldots , a_{in})$ for $2\leq i \leq 2^n$. Also we have at most $2^n$, $G$-layers. First we label the vertices  of one $G$-layer with a distinguishing $D(G)$-labeling. Next we label the $j$th vertex of the $i$th $G$-layer with the $j$th component of the sequence $a_i$ for any $2\leq i \leq m$ and $1\leq j \leq n$. Since we have at most $2^n$, $G$-layers, so if $D(G)=2$ then we can label other $G$-layers such that they have distinct  sequence of labels from the layer  have been labeled distinguishingly.   Thus all automorphism of the strong product $G\boxtimes H$ generated by the automorphism of $G$  are broken, since one of $G$-layers assumes a distinguishing labeling. Also, any two $G$-layers cannot be interchanged as they have different sequence of labels.\qed

%It may happen that $D(G\boxtimes H) < D(G)$ when $|V(H)|\leq 2^{|V(G)|}$. For instance, we showed that $D(C_5\boxtimes P_m)=2$ for any $3\leq m$, while $D(C_5)=3$.

Let $l$ be a positive integer and $A(l,n)$ be the set of sequences $a_i=(a_{i1},\ldots , a_{in})$  such that $1\leq a_{ik}\leq l$ for any $1\leq k \leq n$. Thus $|A(l,n)|=l^n$. By this notation we can obtain an upper bound for strong product of $S$-thin graphs.
\begin{theorem}
Let  $G$ and $H$ be two nonisomorphic simple connected, $S$-thin prime graphs of orders $n$ and $m$, respectively. 
\begin{itemize}
\item[\rm{(i)}] If $D(G)\geq 2$ and $D(G)\neq \lceil {\rm log}_n{m-1} \rceil$, then $D(G\boxtimes H)\leq {\rm  max}\{D(G),\lceil {\rm log}_n{m-1} \rceil\}$. 
\item[\rm{(ii)}] If $D(G)\geq 2$ and $D(G)= \lceil {\rm log}_n{m-1} \rceil$, then  $D(G\boxtimes H)\leq D(G)+1$.
\item[\rm{(iii)}] If $D(G)=1$, then $D(G\boxtimes H)\leq \lceil {\rm log}_n{m} \rceil$.
\end{itemize}
\end{theorem}
\proof (i, ii)  We have $m$, $G$-layers, say $G^{w_1},\ldots , G^{w_m}$, where $V(H)=\{w_1,\ldots , w_m\}$. First we label the vertices  of  $G^{w_1}$-layer with a distinguishing $D(G)$-labeling. Next we label the $j$th vertex of the $i$th $G^{w_i}$-layer with the $j$th component of the sequence $a_i$ for any $2\leq  i \leq  m$ and $1\leq  j \leq  n$.  We use ${\rm max}\{D(G), d\}$ labels where $d ={\rm min}\{l:~|A(l,n)|\geq m-1\}$. With respect to the $|A(l,n)|=l^n$, it can be computed that $d = \lceil {\rm log}_n{m-1} \rceil$.
All automorphism of the strong product $G\boxtimes H$ generated by the automorphism of $G$  are broken, since one of $G$-layers assumes a distinguishing labeling. If $D(G)\neq \lceil {\rm log}_n{m-1} \rceil$, then any two $G$-layers cannot be interchanged as they have different sequence of labels. If $D(G)= \lceil {\rm log}_n{m-1} \rceil$, then we use a new label such that  any two $G$-layers have different sequence of labels and so cannot be interchanged.

(iii) Since $D(G)=1$ so if we label the vertices of each $G^{w_i}$-layer with the the components of the sequence $a_i$, then by a similar argument as (i) we have a distinguishing labeling with $d={\rm min}\{l:~|A(l,n)|\geq m\}$ labels. So we have the result. \qed

\section{Distinguishing index of strong product of two graphs}
In this section we investigate the distinguishing index of strong product of graphs. %Let us start with an easy and useful lemma. 
We say that a graph $G$ is almost spanned by a subgraph $H$ if $G-v$, the graph obtained from $G$ by removal of a vertex $v$ and all edges incident to $v$, is spanned by $H$ for some $v \in  V (G)$. The following two observations will play a crucial role
in this section.
\begin{lemma} {\rm \cite{nord}}\label{nordspanning}
If a graph $G$ is spanned or almost spanned by a subgraph $H$, then $D'(G) \leq D'(H) + 1$.
\end{lemma}
\begin{lemma}\label{distindspann}
Let $G$ be a graph and $H$ be a spanning subgraph of $G$. If ${\rm Aut}(G)$ is a subgroup of ${\rm Aut}(H)$, then $D'(G)\leq  D'(H)$.
\end{lemma}
\proof  If we call the edges of $G$ which are the edges of $H$, $H$-edges, and the others non-$H$-edges, then since ${\rm Aut}(G)\subseteq{\rm Aut}(H)$, we can conclude that each automorphism of $G$ maps $H$-edges to $H$-edges and non-$H$-edges to non-$H$-edges. So assigning  each distinguishing edge labeling of $H$ to $G$ and assigning  non-$H$-edges a repeated label we make    a distinguishing edge labeling of $G$.  \qed

Since for two distinct simple connected, $S$-thin prime graphs we have  ${\rm Aut}(G\boxtimes H) = {\rm Aut}(G\Box H)$, so a direct consequence of Lemmas \ref{nordspanning} and \ref{distindspann} is as follows:
\begin{theorem}\label{disindstrocartesi}
\begin{itemize}
\item[\rm{(i)}] If $G$ and $H$ are two simple connected graphs, then $D'(G\boxtimes H)\leq D'(G\Box H)+1$.
\item[\rm{(ii)}] If $G$ and $H$ are two  simple connected, $S$-thin prime graphs, then $D'(G\boxtimes H)\leq D'(G\Box H)$.
\end{itemize}
\end{theorem}

By the value of the distinguishing index of Cartesian product of $S$-thin graphs in \cite{A. Gorzkowska R. Kalinowski and M. Pilsniak} and Theorem \ref{disindstrocartesi},  we can obtain this value for the strong product of them as the two following corollaries.

\begin{corollary}
\begin{itemize}
\item[\rm{(i)}] The strong product $P_m\boxtimes P_n$ of two paths of orders $m\geq2$ and $n\geq2$ has the distinguishing index equal to two, except $D'(P_2\boxtimes P_2)=3$.
\item[\rm{(ii)}] The strong product $C_m\boxtimes C_n$ of two cycles of orders $m\geq3$ and $n\geq3$ has the distinguishing index equal to two.
\item[\rm{(iii)}] The strong product $P_m\boxtimes C_n$  of orders $m\geq2$ and $n\geq3$ has the distinguishing index equal to two.
\end{itemize}
\end{corollary}
\begin{corollary}
Let $ G$ be a  connected, $S$-thin graph	and $\boxtimes G^k$ the $k$-th power of $G$ with respect to the strong product. Then $D'(\boxtimes G^k) = 2$ for $k \geq 2$.
\end{corollary}

Pil\'sniak in \cite{nord} showed that the distinguishing index of traceable graphs, graphs with a Hamiltonian path, of order equal or  greater than seven is at most two. Also Kr\'al et al. in \cite{Kral etal} showed that if $G_0, \ldots , G_{\Delta -1}$ are non-trivial connected graphs with maximum degree at most $\Delta \geq2$, then $G_0 \boxtimes \cdots \boxtimes G_{\Delta -1}$ is a  Hamiltonian graph. By this note we can prove the following theorem.
\begin{theorem}
Let $G_0, \ldots , G_{\Delta -1}$ be non-trivial connected graphs of orders $n_i$ with maximum degree at most $\Delta \geq2$. If $\prod_{i=0}^{\Delta -1}n_i\geq7$ then $D'(G_0 \boxtimes \cdots \boxtimes G_{\Delta -1})\leq 2$.
\end{theorem}
\begin{corollary}
Let $G$ be a non-trivial connected graph of order $n$ with maximum degree at most $\Delta \geq2$. Then $D'( \boxtimes G^{\Delta})= 2$.
\end{corollary}


\begin{thebibliography}{99}
	
	\bibitem{Albert} M.O. Albertson and K.L. Collins, {\it Symmetry breaking in graphs}, Electron. J. Combin. 3 (1996), \#R18.

	
	\bibitem{Collins and A. N. Trenk} K. L. Collins and A. N. Trenk, {\it The distinguishing chromatic number}, Electron. J. Combin. 13 (2006), \#R16.
	
\bibitem{A. Gorzkowska R. Kalinowski and M. Pilsniak} A. Gorzkowska, R. Kalinowski, and M. Pil\'sniak, {\it The distinguishing index of the Cartesian product of finite
graphs}, Ars Math. Contem. 12 (2017), 77-87.
	
	\bibitem{Sandi} R. Hammack, W. Imrich and S. Kla\v{v}zar, {\it Handbook of product graphs (second edition)}, Taylor \& Francis group 2011.  
	
	
	\bibitem{Imrich1969} W. Imrich, {\it Automorphismen und das kartesische Produkt von Graphen}, O¨ sterreich. Akad.Wiss. Math.-Natur. Kl. S.-B. II 177 (1969), 203-214.
	
	\bibitem{Imrich} W. Imrich, J. Jerebic and S. Kla\v{v}zar, {\it The distinguishing number of Cartesian products of complete graphs}, European J. Combin. 29 (4) (2008), 922-929.
	
	\bibitem{Imrich and  Klavzar} W. Imrich and S. Kla\v{v}zar, {\it Distinguishing Cartesian powers of graphs}, J. Graph Theory, 53.3 (2006), 250-260.

	
	\bibitem{Kali1} R. Kalinowski and M. Pil\'sniak, {\it Distinguishing graphs by edge colourings}, European J. Combin. 45 (2015), 124-131.
	
\bibitem{Kral etal} D. Kr\'al, J.  Maxov\'a,  R. \v{S}\'amal, and P.  Podbrdsk\'y, {\it  Hamilton cycles in strong products of graphs}, J. Graph Theory, 48 (4) (2005), 299-321.
	
	\bibitem{Miller} D. J. Miller, {\it The automorphism group of a product of graphs}, Proc. Amer. Math. Soc. 25 (1970), 24-28.
	

	
	\bibitem{nord}  M. Pil\'sniak, {\it Nordhaus-Gaddum bounds for the distinguishing index}, Available at \texttt{www.ii.uj.edu.pl/preMD/}. 
	

	
\end{thebibliography}
\end{document}